\documentclass[12pt]{amsart}
\usepackage{amssymb}
\usepackage[all]{xypic}

\newtheorem{thm}{Theorem}

\newtheorem{cor}{Corollary}

\DeclareMathOperator{\var}{Var}
\DeclareMathOperator{\covar}{Cov}
\DeclareMathOperator{\p}{pr}

\newcommand{\bH}{\underline{H}}

\renewcommand{\nu}[1]{\overline{#1}}
\newcommand{\un}[1]{\underline{#1}}
\newcommand{\ws}{w_S}

\newcommand{\s}{*=0{\bullet}}
\renewcommand{\c}{*-<1pt>{\circ}}
\newcommand{\n}{{\curlywedge}}
\newcommand{\e}[1]{{\makebox[0pt][r]{\ensuremath{#1\hspace{-0.5ex}}}}}
\newcommand{\fe}{\ar@{-}[luu]}
\newcommand{\ef}{\ar@{-}[ruu]}
\newcommand{\de}{\ar@{--}[luu]}
\newcommand{\ed}{\ar@{--}[ruu]}
\newcommand{\fd}{\fe\ed}
\newcommand{\df}{\ef\de}

\newcommand{\tbl}{\xymatrix@!0@C=3.4ex@R=2.2ex}

\title{Conway's napkin problem}

\author[A. Claesson]{Anders Claesson}
\address{Division of Mathematics,
Department of Chemistry and Biomedical Sciences, University of Kalmar, Sweden}
\email{anders.claesson@hik.se}

\author[T. K. Petersen]{T. Kyle Petersen}
\address{Department of Mathematics, Brandeis University, Waltham, MA, USA, 02454}
\email{tkpeters@brandeis.edu}

\begin{document}
\begin{abstract}
  The napkin problem was first posed by John H. Conway, and written up
  as a `toughie' in ``Mathematical Puzzles: A Connoisseur's Collection,"
  by Peter Winkler. To paraphrase Winkler's book, there is a banquet
  dinner to be served at a mathematics conference. At a particular
  table, $n$ men are to be seated around a circular table. There are $n$
  napkins, exactly one between each of the place settings. Being doubly
  cursed as both men and mathematicians, they are all assumed to be
  ignorant of table etiquette. The men come to sit at the table one at a
  time and in random order. When a guest sits down, he will prefer the
  left napkin with probability $p$ and the right napkin with probability
  $q = 1-p$. If there are napkins on both sides of the place setting, he will
  choose the napkin he prefers. If he finds only one napkin available,
  he will take that napkin (though it may not be the napkin he
  wants). The third possibility is that no napkin is available, and the
  unfortunate guest is faced with the prospect of going through dinner
  without any napkin!

  We think of the question of how many people don't get napkins as a
  statistic for signed permutations, where the permutation gives the
  order in which people sit and the sign tells us whether they initially
  reach left or right. We also keep track of the number of guests who
  get a napkin, but not the napkin they prefer. We find the generating
  function for the joint distribution of these statistics, and use it to
  answer questions like: What is the probability that every guest
  receives a napkin? How many guests do we expect to be without a
  napkin? How many guests are happy with the napkin they receive?
\end{abstract}

\maketitle

\section{Introduction}
The problem studied in this article first appeared in the book
``Mathematical Puzzles: A Connoisseur's Collection," by Peter Winkler
\cite{Winkler}, and was inspired by a true story. Rather than
recounting the problem and the story ourselves, we prefer to quote
directly from ``Mathematical Puzzles":

\begin{quote}
\textsc{The Malicious Maitre D'}

{\small At a mathematics conference banquet, 48 male mathematicians,
none of them knowledgeable about table etiquette, find themselves
assigned to a big circular table. On the table, between each pair of
settings, is a coffee cup containing a cloth napkin. As each person is
seated (by the maitre d'), he takes a napkin from his left or right;
if both napkins are present, he chooses randomly (but the maitre d'
doesn't get to see which one he chose).

In what order should the seats be filled to maximize the expected
number of mathematicians who don't get napkins?

\ldots This problem can be traced to a particular event. Princeton
mathematician John H. Conway came to Bell Labs on March 30, 2001 to
give a ``General Research Colloquium." At lunchtime, [Winkler] found
himself sitting between Conway and computer scientist Rob Pike (now of
Google), and the napkins and coffee cups were as described in the
puzzle. Conway asked how many diners would be without napkins if they
were seated in \emph{random} order, and Pike said: ``Here's an easier
question---what's the \emph{worst} order?"  }
\end{quote}

The problem of the Malicious Maitre D' is not horribly difficult; if
you're having trouble finding a solution, you can see Winkler's book
for a nice explanation. In this paper, it is Conway's problem that we
focus on. Again, from \cite{Winkler}:

\begin{quote}
\textsc{Napkins in a Random Setting}

{\small Remember the conference banquet, where a bunch of
mathematicians find themselves assigned to a big circular table?
Again, on the table, between each pair of settings, is a coffee cup
containing a cloth napkin. As each person sits down, he takes a napkin
from his left or right; if both napkins are present, he chooses
randomly.

This time there is no maitre d'; the seats are occupied in random
order. If the number of mathematicians is large, what fraction of them
(asymptotically) will end up without a napkin?  }
\end{quote}

Let $p$ be the probability that a diner prefers the left napkin and $q=1-p$ be the probability that a diner prefers the right napkin. For the case $p=q=1/2$, Winkler's book gives two proofs of the answer: $(2- \sqrt{e})^2 \approx .12339675$. One is combinatorial, while the other, taken from a more general result due to Aidan Sudbury \cite{Sudbury}, is analytical. In fact Sudbury gives the expected proportion of diners without a napkin as (asymptotically)
\begin{equation}\label{eq:Sud}
\frac{(1 - pe^q)(1 - qe^p)}{pq}.
\end{equation}
(As an aside, we took an informal survey of 55 mathematicians and found about 69\%
would prefer the napkin on the left. According to Sudbury's result, we would
thus expect about 10.58 percent of the guests to get stuck without
napkins.)

In this paper, we use combinatorial methods to produce the generating
function for the probability that at a table for $n$ people, $i$ of
them have no napkin and $j$ of them have a napkin, but not the napkin
they prefer. This generating function allows for a thorough
statistical analysis of the problem, including a new proof of \eqref{eq:Sud}.

>From our point of view, the number of people without a napkin is a
statistic for signed permutations; just not one so well studied as
inversions, descents, and such. We consider the order in which guests
sit down at a place as a permutation of $[n]= \{1,2,\ldots,n\}$, while
their preference for the right or left napkin is given by a plus or
minus sign. We label the places $1,2,3,\ldots,n$ counter-clockwise, so
that place $i$ has place $i-1$ on its left, place $i+1$ on its right,
and place $n$ is to the left of place 1. With this convention, the signed permutation $(2,-3,4,-1)$ describes the following sequence of events at a table for four. The person sitting in place 4 sits down
first and takes the napkin on his left. The person in place 1 sits
next and takes the napkin on his right. The person at
place 2 sits third and wants to take the napkin on his left, but since that napkin is already taken, he is forced to take the napkin on the right. Finally, poor person 3 sits last to find no napkin at his place.

\begin{figure}\label{fig:intro_ex}
$$
\begin{array}{|l|l|}
\hline
\xymatrix@!0@C=4ex@R=2.9ex{
  \\
  &\n && \n && \n && \n &&&\\
  \\
  && \s    && \s  && \s    && \s\fe \\
  && \e{2} && \e{-3} && \e{4} && \e{-1} \\
  \\
}
&
\xymatrix@!0@C=4ex@R=2.9ex{
  \\
  &\n && \n && \n && \n &&&\\
  \\
  && \s\ef && \s  && \s    && \s\fe \\
  && \e{2} && \e{-3} && \e{4} && \e{-1} \\
  \\
}\\
\hline
\xymatrix@!0@C=4ex@R=2.9ex{
  \\
  &\n && \n && \n && \n &&&\\
  \\
  && \s\ef && \s\df  && \s && \s\fe \\
  && \e{2} && \e{-3} && \e{4} && \e{-1} \\
  \\
}&
\xymatrix@!0@C=4ex@R=2.9ex{
  \\
  &\n && \n && \n && \n &&&\\
  \\
  && \s\ef && \s\df  && \s\ed && \s\fe \\
  && \e{2} && \e{-3} && \e{4} && \e{-1} \\
  \\
}\\
\hline
\end{array}
$$
\caption{A table for four.}
\end{figure}
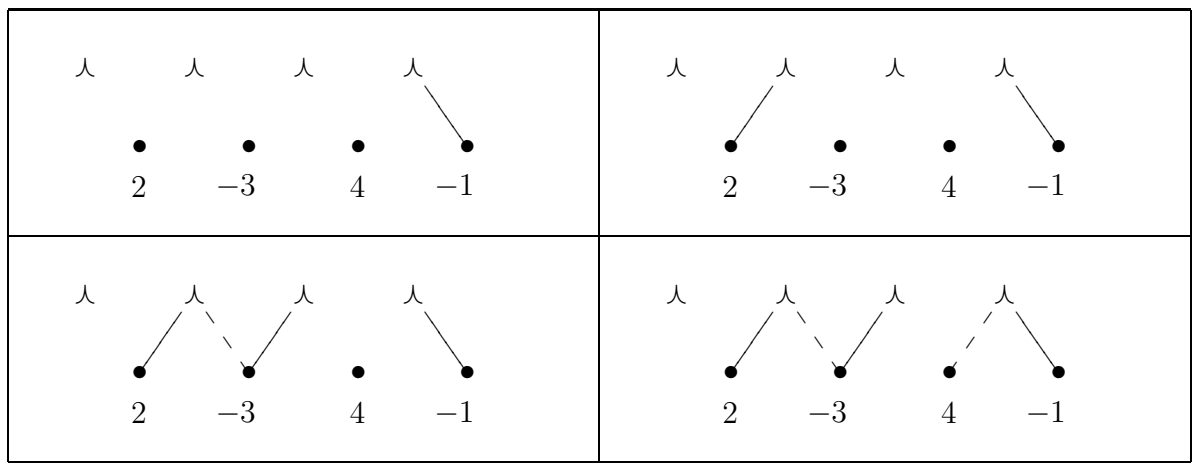

In Figure~\ref{fig:intro_ex} we illustrate this sequence of events. There
a ``$\bullet$'' and a ``$\curlywedge$'' represent a seat and a
napkin, respectively. Moreover, a solid line symbolizes an arm taking
a napkin and a broken line symbolizes an imaginary arm reaching for
the preferred napkin. When looking at these diagrams, keep in mind that
the table is circular: the guests at places $1$ and $4$ are neighbors.
Each signed permutation thus corresponds to a particular set of diners
who will have no napkin. In particular, for each permutation there is
some number of diners without napkins, and we seek to determine what
proportion of signed permutations leave a given number of diners
napkinless.

We remark that this problem makes for interesting mealtime
conversation, and we have heard many suggestions for variations on the
theme. What if instead of a math conference banquet with all men, it
is a dinner party for couples, and the couples enter one at a time,
lady sitting first? What if the dinner party is a mixer for singles,
where now all the ladies enter first, sitting in alternating
seats?\footnote{This variation is actually not very interesting
mathematically, as it is not difficult to see the women will always
get their choice of napkin, and that 1/4 (or $pq$ in general) of the men are expected to be
without a napkin.} What if the guests don't mind looking farther
afield than simply to their immediate left and right, say, reaching as
many as two places over for a napkin? Exploring answers to these
questions may lead to some interesting mathematics, but in this
article we stay within the original question, where the guests are all
male mathematicians, they enter in a totally random order, and they
are too shy to reach beyond their immediate left or right to find a
napkin.

In section \ref{sec:def} we define the generating functions that we
use throughout the paper. Section \ref{sec:every} explores the
question of when everyone gets a napkin, and makes connections with
combinatorial objects called \emph{ordered bipartitions}, due to
Dominique Foata and Doron Zeilberger \cite{FoataZeil}. This connection
makes subsequent proofs much simpler, as in section
\ref{sec:ordgenfun}, where we derive powerful identities involving our
generating functions that ultimately lead to exact formulas. In
section \ref{sec:expect} we answer the original question of how many
guests are expected to be without a napkin, as well as provide some
other statistics of interest. Section \ref{sec:alternate} gives a
second derivation of the generating function for the expected number
of guests without a napkin.

\section{Definitions}\label{sec:def}

Let $p$ be the probability that a diner will reach to the left, and
$q = 1-p$ be the probability of reaching right. We denote the set of signed permutations of $[n]$ by
$\mathcal{C}_n$. For any $\pi \in \mathcal{C}_n$, let $|\pi|_{-}$ and
$|\pi|_{+}$ be the number of negative and positive entries in $\pi$,
respectively. Also, let $|\pi|=|\pi|_{-} + |\pi|_{+} = n$ be the length
of $\pi$. Let $o(\pi)$ be the number of guests without a napkin after every guest has been seated as
described by $\pi$. Note that the number of guests without a napkin is
equal to the number of napkins left on the table. Furthermore, let
$m(\pi)$ be the number of people who get a napkin, but not their first
choice. We say a guest is \emph{napkinless} if he has no napkin,
and a guest is \emph{frustrated} if he gets a napkin, but not the
napkin he originally wanted.\footnote{When discussing this problem,
French mathematician Sylvie Corteel argued that if the mathematicians
were French, they would never take the ``incorrect" napkin. If the
napkin they wanted was not there, they would simply cross their arms
and refuse to eat. But of course, that is a different problem (how
does the number of napkinless guests change as the proportion of
French diners changes?).} Otherwise, we say the guest is
\emph{happy}. Define the \emph{weight} of $\pi$ as
$$ w(\pi) := p^{|\pi|_{-}}q^{|\pi|_{+}}x^{o(\pi)}y^{m(\pi)}.
$$
For example, if $\pi = (2, -1, 3, 4)$, then $|\pi|_{-} = 1$,
$|\pi|_{+} = 3$, and $|\pi| = 4$. Also, the napkin between places 2 and 3 is unused
(guest 4 is the unlucky one), and although person 1 gets a napkin, it
was not the one he wanted. Thus, $o(\pi) = 1$, $m(\pi) = 1$, and two
of the guests are happy: $w(\pi) = pq^3xy$.

The generating function we are interested in is
\begin{equation}\label{eq:Cgf}
  C(p; x,y,z)
  := \sum_{n \geq 1}\sum_{\pi \in \mathcal{C}_n}w(\pi)\frac{z^n}{n!} = \sum_{\substack{ i,j\geq0 \\ n \geq 1}} \p(i,j,n)x^i y^j z^n,
\end{equation}
where $\p(i,j,n)$ denotes the probability that at a table for $n$ people,
$i$ of them are napkinless, and $j$ of them are frustrated.

The main result of this paper, Theorem \ref{thm:exact}, is an exact formula for
$C(p; x,y,z)$. Our approach is to first ``straighten" the table.

Suppose that instead of a circular table with $n$ places and $n$
napkins, we look at a straight table with $n$ places and $n+1$
napkins, so that each place has a napkin on both its left and
right. If we know all that can possibly happen in this situation, then
in order to determine the circular case, we just consider that the
last person to enter the room sits ``between" the first and last
person on the linear table. Let us make this connection more precise.

Let $\mathcal{N}_n$ be the set of all signed permutations of $[n]$
that result in taking neither napkins from the ends of the
table. Similarly, $\mathcal{L}_n$ (resp. $\mathcal{R}_n$) denotes the
set of all signed permutations of $[n]$ that result in the left end
napkin being taken but not the right (resp. right but not left), and
$\mathcal{B}_n$ denotes those signed permutations that result in both
end napkins being taken.  We note that $\mathcal{C}_n = \mathcal{N}_n
\cup \mathcal{L}_n \cup \mathcal{R}_n \cup \mathcal{B}_n$.

For the straight table we define the weight of $\pi$ as
$$ \ws(\pi) := p^{|\pi|_{-}}q^{|\pi|_{+}}x^{o(\pi)}y^{m(\pi)},
$$ where $o(\pi)$ is still the number of people without a napkin at
the (now linear) table, so that $o(\pi)+1$ is the number of napkins
left on the table. Moreover, we define the following generating
functions for the linear table:
\begin{align*}
  S(p;x,y,z) &:= \sum_{n\geq 0} \sum_{\pi \in \mathcal{C}_n} \ws(\pi) \frac{z^{n}}{n!};\\
  N(p;x,y,z) &:= \sum_{n\geq 0} \sum_{\pi \in \mathcal{N}_n} \ws(\pi) \frac{z^{n}}{n!};\\
  L(p;x,y,z) &:= \sum_{n\geq 1} \sum_{\pi \in \mathcal{L}_n} \ws(\pi) \frac{z^{n}}{n!};\\
  R(p;x,y,z) &:= \sum_{n\geq 1} \sum_{\pi \in \mathcal{R}_n} \ws(\pi) \frac{z^{n}}{n!};\\
  B(p;x,y,z) &:= \sum_{n\geq 1} \sum_{\pi \in \mathcal{B}_n} \ws(\pi) \frac{z^{n}}{n!}.
\end{align*}
By construction, we have
\begin{equation}\label{eq:SNLRB}
  S(p;x,y,z) = N(p;x,y,z) + L(p;x,y,z) + R(p;x,y,z) + B(p;x,y,z).
\end{equation}
Further, we can observe that by symmetry
\begin{alignat*}{3}
  C(p;x,y,z) &= C(q;x,y,z), &\;\,&\;\,& S(p;x,y,z) &= S(q;x,y,z), \\
  N(p;x,y,z) &= N(q;x,y,z), &\;\,&\;\,& B(p;x,y,z) &= B(q;x,y,z),
\end{alignat*}
and
\begin{equation}
  R(p;x,y,z) = L(q;x,y,z).\label{eq:RL}
\end{equation}
>From now on we will usually suppress the $p$ and write $C(x,y,z)$
instead of $C(p;x,y,z)$, $S(x,y,z)$ instead of $S(p;x,y,z)$, etc.

Now let us recast the generating function $C(x,y,z)$ in terms of the
generating functions for the linear table. Everything that has
happened before the last person sits down at a table for $n$ people
can be considered the result of a signed permutation of $[n-1]$
playing out on a linear table. If the last person sits at place $n$, this is obvious, but in general the labeling of the seats is unimportant, i.e., we can always cyclically permute the labels on the seats without changing the weight of the permutation.

In what follows, $\pi = \sigma\tau$ means that $\pi$ is the concatenation of
$\sigma$ with $\tau$. For instance, $(3,-1)(-6)(2,-4,5) = (3,-1,-6,2,-4,5)$.

\begin{figure}
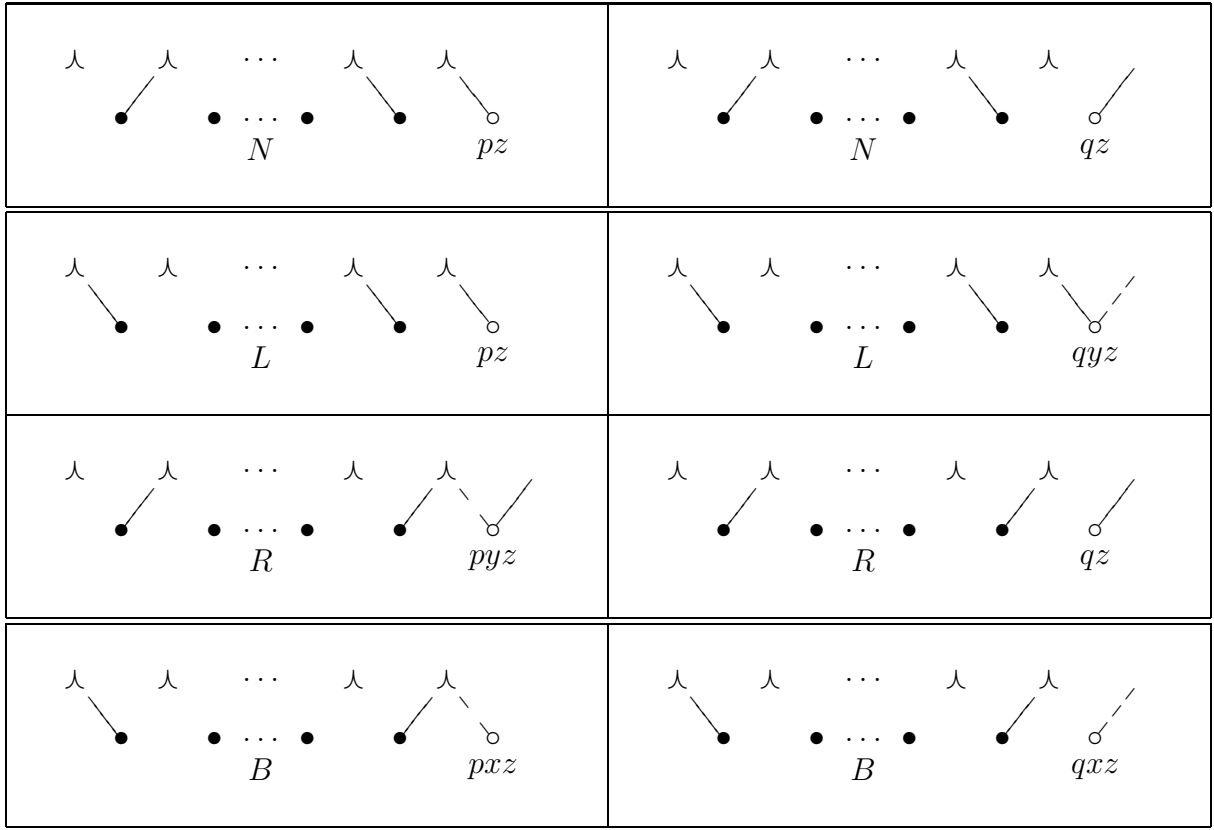
\label{fig:C_NLRB}
$$
\begin{array}{|c|c|}
\hline
&\\
\tbl{
  &\n && \n && \ldots &&  \n && \n &&& \\
    \\
  && \s\ef && \s & \ldots & \s && \s\fe && \c\fe \\
  &&&&& N &&&&& pz \\
}&
\tbl{
  &\n && \n && \ldots &&  \n && \n &&& \\
    \\
  && \s\ef && \s & \ldots & \s && \s\fe && \c\ef \\
  &&&&& N &&&&& qz \\
}\\
&\\ \hline\hline &\\
\tbl{
  &\n && \n && \ldots &&  \n && \n &&& \\
    \\
  && \s\fe && \s & \ldots & \s && \s\fe && \c\fe \\
  &&&&& L &&&&& pz \\
}&
\tbl{
  &\n && \n && \ldots &&  \n && \n &&& \\
    \\
  && \s\fe && \s & \ldots & \s && \s\fe && \c\fd \\
  &&&&& L &&&&& qyz \\
}\\
&\\ \hline &\\
\tbl{
  &\n && \n && \ldots &&  \n && \n &&& \\
    \\
  && \s\ef && \s & \ldots & \s && \s\ef && \c\df \\
  &&&&& R &&&&& pyz \\
}&
\tbl{
  &\n && \n && \ldots &&  \n && \n &&& \\
    \\
  && \s\ef && \s & \ldots & \s && \s\ef && \c\ef \\
  &&&&& R &&&&& qz \\
}\\
&\\ \hline\hline &\\
\tbl{
  &\n && \n && \ldots &&  \n && \n &&& \\
    \\
  && \s\fe && \s & \ldots & \s && \s\ef && \c\de \\
  &&&&& B &&&&& pxz \\
}&
\tbl{
  &\n && \n && \ldots &&  \n && \n &&& \\
    \\
  && \s\fe && \s & \ldots & \s && \s\ef && \c\ed \\
  &&&&& B &&&&& qxz \\
}\\
& \\ \hline
\end{array}
$$
\caption{$C = z\big( N + (p+qy)L + (py+q)R + xB \big)$.}
\end{figure}

Let us consider a circular table for $n$ guests and assume that $\pi =
\sigma(\pm n)\tau$ is a member of $\mathcal{C}_n$. Let $\pi' =
\tau\sigma$. If the last person walks in and has both napkins
available, then $\pi'$ must have resulted in leaving both end napkins
on a linear table, that is, $\pi' \in \mathcal{N}_{n-1}$. Whether the
last guest prefers the left napkin or the right napkin, he will get
his choice. So, the weight of $\pi$ is:
$$
w(\pi) =
\begin{cases}
p\ws(\pi') & \mbox{if person } n \mbox{ prefers left,}\\
q\ws(\pi') & \mbox{if person } n \mbox{ prefers right.}
\end{cases}
$$
If person $n$ walks in to find only the left napkin available,
$\pi' \in \mathcal{L}_{n-1}$, and the last person
will take that napkin regardless of preference, getting the one he
wants with probability $p$, and getting frustrated with probability $q$:
$$w(\pi) =
\begin{cases}
p\ws(\pi') & \mbox{if person } n \mbox{ prefers left,}\\
qy\ws(\pi') & \mbox{if person } n \mbox{ prefers right.}
\end{cases}
$$
Similarly, if person $n$ walks in to find only the right napkin available then
$$w(\pi) =
\begin{cases}
py\ws(\pi') & \mbox{if person } n \mbox{ prefers left,}\\
q\ws(\pi') & \mbox{if person } n \mbox{ prefers right.}
\end{cases}
$$
Finally, guest $n$ can walk in to find both
napkins already taken. In this case, we know that $\pi' \in
\mathcal{B}_{n-1}$ and that guest $n$ will be one of the napkinless
guests:
$$w(\pi) =
\begin{cases}
px\ws(\pi') & \mbox{if person } n \mbox{ prefers left,}\\
qx\ws(\pi') & \mbox{if person } n \mbox{ prefers right.}
\end{cases}
$$
Figure~\ref{fig:C_NLRB} illustrates these, altogether 8,
possibilities in the special case when the last person to sits at
place $n$. All other cases can be reduced to this case by cyclically
shifting the picture.

In terms of generating functions we have showed that
\begin{equation}\label{eq:CNLRB}
  C(x,y,z) = z\big( N(x,y,z) + (p+qy)L(x,y,z) + (py+q)R(x,y,z) + xB(x,y,z) \big).
\end{equation}


\section{A warm-up problem: when does everyone have a napkin?}
\label{sec:every}

A natural question to ask (perhaps easier than the general question)
is: what is the probability, $\p(0,1,n) = \p(0,n)$, that at a circular table for
$n$ people, every guest has a napkin? The generating function for this probability
is $C(0,y,z)$. By equation \eqref{eq:CNLRB}, we have
\[C(0,y,z) = z\big( N(0,y,z) + (p+qy)L(0,y,z) + (py+q)R(0,y,z) \big),
\]
but since (on a straight table with at least one person) there is
always at least one person without a napkin if neither end napkin is
taken, $N(0,y,z) = 1$, and the above equation reduces to
\begin{equation}\label{eq:C0}
  C(0,y,z) = z\big( 1 + (p+qy)L(0,y,z) + (py+q)R(0,y,z) \big).
\end{equation}
This makes intuitive sense because if everyone is to have a napkin,
they all need to take the left napkin, or all need to take the right
napkin. Therefore we turn our attention to $L(0,y,z)$ (since $R(0,y,z)$ follows by swapping $p$ and $q$).

For $n$ a positive integer, let us consider a
straight table for $n$ guests and let $\pi\in\mathcal{L}_n$. Further,
let $\pi = \sigma(\pm n)\tau$, where $\sigma$ and $\tau$ are the
signed permutations consisting of the letters to the left and to the
right of $n$, respectively.

If the last person to enter sits down at the rightmost seat (i.e.,
$\tau$ is empty) then, since $\pi\in\mathcal{L}_n$, he necessarily
prefers the left napkin. So,
$$ \ws(\pi) = \ws(\sigma)p.
$$

Assume that the last person to enter sits down with at least one
guest to his right (i.e., $\tau$ is not empty). Then he will be
happy with probability $p$ and frustrated with probability $q$:
$$ \ws(\pi) =
\begin{cases}
\ws(\sigma)p\ws(\tau) & \mbox{if person } n \mbox{ prefers left,}\\
\ws(\sigma)qy\ws(\tau) & \mbox{if person } n \mbox{ prefers right.}
\end{cases}
$$

\begin{figure}
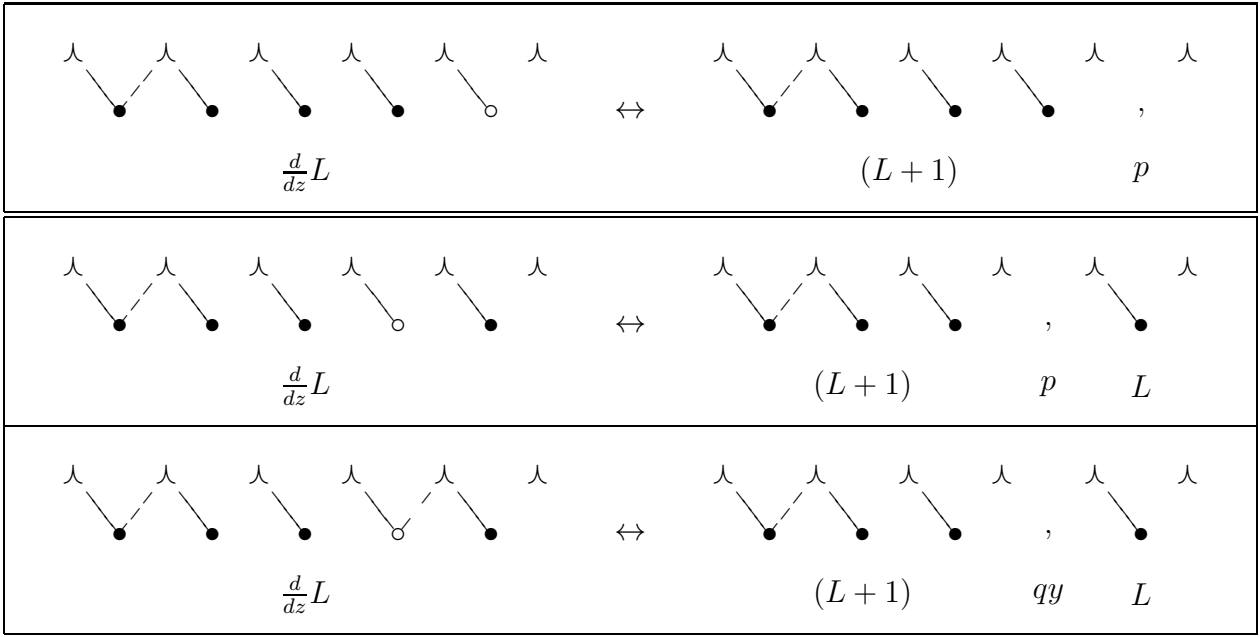
\label{fig:L0}
$$
\begin{array}{|c|}
\hline
\tbl{
  \\
  &\n && \n && \n &&  \n && \n && \n &&&  &\n && \n && \n && \n &&\n && \n &\\
  \\
  && \s\fd && \s\fe && \s\fe && \s\fe && \c\fe &&& \leftrightarrow &&& \s\fd && \s\fe && \s\fe &&\s\fe&&,&&  \\
  \\
  &&&&&& \frac{d}{dz}L &&&&&&&&&&&&&(L+1)&&&&& p\\
  \\
}\\
\hline
\hline
\tbl{
  \\
  &\n && \n && \n &&  \n && \n && \n &&&  &\n && \n && \n && \n &&\n && \n &\\
  \\
  && \s\fd && \s\fe && \s\fe && \c\fe && \s\fe &&& \leftrightarrow &&& \s\fd && \s\fe && \s\fe &&,&& \s\fe &&  \\
  \\
  &&&&&& \frac{d}{dz}L &&&&&&&&&&&&(L+1) &&&& p && L\\
  \\
}\\
\hline
\tbl{
  \\
  &\n && \n && \n &&  \n && \n && \n & && &\n && \n && \n && \n &&\n && \n &\\
  \\
  && \s\fd && \s\fe && \s\fe && \c\fd && \s\fe &&& \leftrightarrow &&& \s\fd && \s\fe && \s\fe &&,&& \s\fe &&  \\
  \\
  &&&&&& \frac{d}{dz}L &&&&&&&&&&&& (L+1) &&&& qy && L\\
  \\
}\\
\hline
\end{array}
$$
\caption{ $\frac{d}{dz}L = (L+1)p + (L+1)pL + (L+1)qyL$, where $L=L(0,y,z)$.}
\end{figure}

Thus the function $L(0,y,z)$ satisfies the following differential equation:
\begin{equation}\label{eq:L0deq}
  \frac{d}{dz}L(0,y,z)
  = \big( L(0,y,z) + 1 \big)\big( p + (p+qy)L(0,y,z) \big).
\end{equation}
Se Figure~\ref{fig:L0} for an illustration.

Let $L_n = [z^n]L(1/2;0,1,z)$, where $[z^n]F(z)$ denotes the
coefficient of $z^n$ in $F(z)$. Using \eqref{eq:L0deq} we compute the
$L_n$ for some small values: $L_1 = 1$, $L_2 = 3$, $L_3 = 13$, $L_4 =
75$, $L_5 = 541$, \dots, and plug them into Sloane's Encyclopedia
\cite{Sloane}. Luckily, we get a hit with sequence A000670! These
numbers happen to be fairly well known as the ``ordered Bell numbers,"
or the number of ordered set partitions of $[n]$. We shall now give a
bijection between the ordered set partitions of $[n]$ and permutations
for which everyone takes the left napkin. This bijection will also
give the general solution to \eqref{eq:L0deq}.

\subsection{Ordered set partitions}

We now give a bijection between ordered set
partitions of $[n]$ and signed permutations that correspond to
everyone at a linear table taking the napkin on the left. An \emph{ordered
set partition}, $\alpha$, of $[n]$ is a word
$$\alpha = B_1 B_2 \cdots B_k
$$
where the ``blocks" $B_1$, $B_2$, \dots, $B_k$ are subsets of $[n]$ such that
$\{B_1,B_2,\ldots,B_k\}$ is a set partition of $[n]$. By convention we
require that the elements of each block are written in decreasing order.

We describe
the bijection with an example. Let
\[\alpha = \{5,2\}\{6\}\{7,4,1\}\{3\}.
\]
First, we give a minus sign to the least element of each block, then
we remove the braces, to obtain
\[\pi = (5,-2,-6,7,4,-1,-3),
\]
a permutation corresponding to a situation where everyone takes the
napkin on the left. We can see that we indeed have produced a
permutation in which everyone takes the napkin on the left by thinking
of $\pi$ as a diagram for the entry times and napkin preferences of
everyone at the table. With $\pi$ as above, we see that the person
seated second from the right entered first and took the napkin on the
left. The person seated third from the right wanted to take that
napkin, but entered fourth, and so was forced to take the napkin on
his left. Because we required that the elements of a block are in
decreasing order, anytime someone has a preference for the napkin on
the right, they find that it was taken by the person to their right.

Clearly this process is reversible. Given a permutation where everyone
takes the napkin on the left, it must have at least one minus sign. In
particular, it must have a minus sign on 1, since the first person to
enter must take the napkin on his left. Anybody with a plus sign
immediately to the left of a person with a minus sign must enter after
that person. And if two or more people with plus signs are sitting
(consecutively) to the left of a person with a minus sign, they must
arrive in order of closeness to the minus sign; the closest to the
minus sign first, followed by the second closest to the minus sign,
and so on. This gives us the block structure of the partition on
$[n]$. Reading $\pi$ from left to right, we separate the blocks by
putting walls immediately to the right of any number appearing with a
minus sign. For example, if
\[\pi = (7, -5, -6, 4, -1, 3, -2),
\]
we convert this permutation into
\[\alpha = \{7,5\}\{6\}\{4,1\}\{3,2\}.
\]

Now we can build the generating function for $L(0,y,z)$ by purely
combinatorial means. Taking an approach from \cite{FoataZeil} (more on
that paper in a bit), let
$$ H(p; y,z)
:= \frac{pz}{1!} + \frac{pqyz^2}{2!} + \cdots + \frac{p(qy)^{n-1}z^n}{n!} + \cdots
= \frac{p(e^{qyz}-1)}{qy}
$$
Then $H(y,z) := H(p; y,z)$ is the generating function for single blocks, $\{n, \ldots, 2,
1\}$, since for every block, only the person corresponding to the
least number gets the napkin he wants, leaving the other $n-1$ people
frustrated. For example with $n = 5$:
$$ w\Bigg(
\vcenter{\tbl{
  \n && \n && \n && \n && \n \\
    \\
  & \s\fd && \s\fd && \s\fd && \s\fd && \s\fe \\
  & qy    && qy    && qy    && qy    && p
}}\,
\Bigg) = p (qy)^4.
$$
Therefore
$$(1-H(y,z))^{-1} = 1 + H(y,z) + (H(y,z))^2 + \cdots
$$
is the generating function for ordered sequences of blocks, and thus
\[ L(0,y,z) = (1 - H(y,z))^{-1} - 1 =
\frac{p(e^{qyz}-1)}{qy-p(e^{qyz}-1)}.
\]
One can easily check that this expression indeed satisfies the
differential equation~\eqref{eq:L0deq}. Due to equations \eqref{eq:RL} and
\eqref{eq:C0} we are now in a position to give the generating function
for the probability that everyone gets a napkin on a circular table:
\begin{align*}
C(0,y,z)
&= z\big( 1 + (p+qy)L(p;0,y,z) + (py+q)L(q;0,y,z) \big) \\
&= z\!\left(1 + \frac{(p+qy)p(e^{qyz}-1)}{qy-p(e^{qyz}-1)}
+ \frac{(py+q)q(e^{pyz}-1)}{py-q(e^{pyz}-1)} \right).
\end{align*}

\subsection{Ordered bipartitions}

Now we will generalize the correspondence just described. The paper of Foata and Zeilberger \cite{FoataZeil} introduces objects called \emph{ordered bipartitions}, which are easiest to think of as
ordered set partitions with some of the subsets underlined. A
\emph{compatible} bipartition is an ordered bipartition where all the
underlined subsets are on the right, e.g.,
\[\alpha =\{5,2\}\{6\}\underline{\{1,4,7\}}\, \underline{\{3\}},
\]
where we adopt the convention that underlined subsets have their
elements written in ascending order. The bijection works as
follows. For every non-underlined group in $\alpha$, we perform the same
operation as above, while for every underlined group we perform the
``opposite" operation. Specifically, we put minus signs on all
\emph{but} the least element before removing the braces. This produces
\[ \pi = ( 5, -2, -6, 1, -4, -7, 3 ),
\]
a permutation where everyone on a linear table receives a napkin. All
we really need to observe is that the underlined groups correspond to
the part of the table where people all take napkins on their right,
while the non-underlined groups all take napkins on the left. Thus,
\begin{align}
  S(0,y,z) &= (L(0,y,z)+1)(R(0,y,z)+1)\nonumber \\
  &= \frac{pqy^2}{\big(q(e^{pyz}-1)-py\big)\big(p(e^{qyz}-1)-qy\big)} \label{eq:S0}
\end{align}

We can take this more general correspondence and see now that the
permutations for which everyone takes the napkin on the left
correspond to the ordered bipartitions with no underlined subsets, the
ordered bipartitions with all subsets underlined correspond to
permutations where everyone takes the right napkin, and the compatible
bipartitions with at least one underlined subset and one
non-underlined subset correspond to the permutations where everyone
gets a napkin and both end napkins are taken. This last observation
gives
\[ B(0,y,z) = L(0,y,z)R(0,y,z).
\]

\section{Ordered bipartitions and generating functions}\label{sec:ordgenfun}

Using ordered bipartitions (not simply the compatible ones), we can encode more than just those permutations where everyone gets a napkin. Using the algorithm below, we can encode the set of \emph{all} signed permutations. Let $\varphi_S$ be this map, where the subscript $S$ reminds us that we are dealing with the
straight table. The image of the injection $\varphi_S$ will lead us to
the main theorem of this section, which gives some wonderful
relationships between the generating functions for the linear
table. Given a signed permutation $\pi$ of $[n]$, we form its image
$\varphi_S(\pi)$ as follows:

\begin{itemize}
\item[(1)] Find the least element, $\pi(i)$ (ignoring signs), that is
  not already included in some subset.

\item[(2a)] If $\pi(i)$ is positive, then underline the set including
  $\pi(i),$ and set $j = i+1$.\\
  While $|\pi(j)| > |\pi(j-1)|$ and $\pi(j)$ negative,\\
  $\mbox{\quad}$ add $\pi(j)$ to the set containing $\pi(i)$, and set $j = j+1$.

\item[(2b)] If $\pi(i)$ is negative, then set $j = i-1$.\\
  While $|\pi(j)| > |\pi(j+1)|$ and $\pi(j)$ positive,\\
  $\mbox{\quad}$ add $\pi(j)$ to the set containing $\pi(i)$, and set $j = j-1$.

\item[(3)] If every element is contained in a set, then delete all
  minus signs and quit. Else, go to (1).
\end{itemize}

Clearly, no two permutations can be mapped to the same bipartition. We
will demonstrate the injection with an example. Start with \[\pi = (9,
1, -3, 2, 5, 6, -4, -7, 8).\] The number 1 is the least element, and
it is positive, so its set will be underlined. We start searching to
the right of 1, looking for negative numbers with absolute value bigger than
1. We get
\[
9 \;\; \underline{\{\,1, -3\,\}} \;\; 2 \;\; 5 \;\; 6 \;\; {-4} \;\; {-7}
\;\; 8
\]
Now 2 is the least element, and it is also positive,
but there are no negative numbers immediately to its right, so it
forms a singleton set,
\[
9 \;\; \underline{\{\,1, -3\,\}} \; \underline{\{\,2\,\}} \;\; 5
\;\; 6 \;\; {-4} \;\; {-7} \;\; 8
\]
Now 4 is the least element not already in a set, and it is
negative. So we start searching to the left of 4, looking for positive
numbers with bigger absolute value. We get
\[
9 \;\; \underline{\{\,1, -3\,\}} \; \underline{\{\,2\,\}} \;\; 5 \;\; \{\,6 ,
-4\,\} \;\; {-7} \;\; 8
\]
The next steps give,
\[
9 \;\; \underline{\{\,1, -3\,\}} \; \underline{\{\,2\,\}} \;
\underline{\{\,5\,\}} \; \{\,6 , -4\,\} \;\; {-7} \;\; 8
\]
\[
9 \; \underline{\{\,1, -3\,\}} \; \underline{\{\,2\,\}} \;
\underline{\{\,5\,\}} \; \{\,6 , -4\,\} \; \{\,-7\,\} \;\; 8
\]
\[
9 \;\; \underline{\{\,1, -3\,\}} \; \underline{\{\,2\,\}} \;
\underline{\{\,5\,\}} \; \{\,6 , -4\,\} \; \{\,-7\,\} \; \underline{\{\,8\,\}}
\]
\[
\underline{\{\,9\,\}} \; \underline{\{\,1, -3\,\}} \; \underline{\{\,2\,\}}
\; \underline{\{\,5\,\}} \; \{\,6 , -4\,\} \; \{\,-7\,\}
\;\underline{\{\,8\,\}}
\]
Then we drop the minus signs to get,
\[
\underline{\{\,9\,\}} \; \underline{\{\,1, 3\,\}} \; \underline{\{\,2\,\}}
\; \underline{\{\,5\,\}} \; \{\,6 , 4\,\} \; \{\,7\,\} \;
\underline{\{\,8\,\}}
\]
so that we have an ordered bipartition of $[n]$ where we write
underlined sets in increasing order.

We will now explain how the weight $w(\pi)$ of the signed permutation
$\pi$ can be read from the corresponding ordered bipartition $\alpha =
\varphi_S(\pi)$. To start with, the number of napkinless diners, $o(\pi)$,
is exactly the number of occurrences of an underlined set followed
immediately by a non-underlined set. The frustrated diners are those
who are not the least element in a block, less the people without any
napkin:
$$m(\pi) = |\alpha| - \ell(\alpha) - o(\pi),
$$ where $|\alpha|$ is the number of elements of the underlying set
of $\alpha$ (same as $|\pi|$ here), and $\ell(\alpha)$ is the number of
blocks in $\alpha$. Let $\un{\alpha}$ (resp. $\nu{\alpha}$) be the
the bipartition formed from the blocks of $\alpha$ that are underlined
(resp. non-underlined). For underlined blocks of $\alpha$ we have that
the least element is positive in $\pi$ and the other elements are negative
in $\pi$. Similarly, for non-underlined blocks of $\alpha$,
the least element is negative in $\pi$ and the other elements are positive
in $\pi$. Thus,
\begin{align*}
  |\pi|_{-} &= |\un{\alpha}| - \ell(\un{\alpha}) + \ell(\nu{\alpha}); \\
  |\pi|_{+} &= |\nu{\alpha}| - \ell(\nu{\alpha}) + \ell(\un{\alpha}). \\
\end{align*}

Let $\varphi_S(\mathcal{C}_n)$ be the image set of
bipartitions corresponding to all signed permutations of $[n]$. By
examining the algorithm describing $\varphi_S$, we see that the set
$\varphi_S(\mathcal{C}_n)$ consists of all ordered bipartitions that never
contain the patterns
$$
\underline{\cdots a\}} \, \{b\} \cdots
\quad\text{or}\quad
\cdots \underline{\{b\}} \, \{ a \cdots, \qquad\text{where $a < b$}.
$$

To simplify notation in what follows, we write $C$ for $C(x,y,z)$, $S$
for $S(x,y,z)$, etc. Further, we write $H$ for $H(p;y,z)$ and $\bH$
for $H(q;y,z)$. The following theorem takes advantage of our combinatorial model. Its power is reflected in the subsequent results it implies.

\begin{thm}\label{thm:BL}
  We have the following formulas:
  \begin{align}
    L + B &= H S \\
    B &= H S \bH
  \end{align}
\end{thm}

\begin{proof}
  If we add a non-underlined block of size $r$ to the left of any
  bipartition in $\varphi_S(\mathcal{C}_n)$, then clearly we get a
  bipartition in $\varphi_S(\mathcal{C}_{n+r})$ that corresponds to a
  permutation in $\mathcal{L}_{n+r} \cup
  \mathcal{B}_{n+r}$. Furthermore, this new block will not change the
  number of people without a napkin (occurrences of underlined blocks
  immediately to the left of non-underlined blocks), and the number of
  new people who get a napkin they don't want is exactly
  $r-1$. Therefore,
  \[ L + B = H S.
  \]

  Similarly, if we add a non-underlined block of size $r$ to the left of
  a a bipartition in $\varphi_S(\mathcal{C}_n)$, and an underlined block of size
  $s$ to the right, then we get a bipartition in $\varphi_S(\mathcal{C}_{n+r+s})$
  that corresponds to a permutation in $\mathcal{B}_{n+r+s}$. Thus,
  \[B = H S \bH,
  \]
  which concludes the proof.
\end{proof}

\begin{cor}\label{thm:KBLNJS}
  We have the following formulas:
  \begin{align*}
    N &= S \big(1-H\big) \big(1-\bH\big)\\
    L &= S H \big(1-\bH\big)\\
    R &= S \big(1-H\big) \bH\\
    B &= S H \bH \\
    C &= z S
    \big( 1 + q(y-1)H + p(y-1)\bH + (x-y)H \bH \big)
  \end{align*}
\end{cor}

\begin{proof}
  From Theorem~\ref{thm:BL} we immediately get the formulas for $B$
  and $L$. Plug these formulas into \eqref{eq:RL} and
  \eqref{eq:SNLRB} and the formulas for $R$ and $N$ follows. Finally,
  \eqref{eq:CNLRB} yields the formula for $C$.
\end{proof}

What Corollary \ref{thm:KBLNJS} tells us is that if we can find an
explicit formula for $S$, then we will have explicit formulas
for all the other generating functions. We will derive such an
explicit formula shortly. First we need to introduce the notion of a
\emph{cyclic} bipartition. A cyclic bipartition is a bipartition for
which one element is distinguished, and only the cyclic ordering of
the blocks matters. As a convention, we put the block containing the
distinguished element at the far right if that block is not
underlined, or at the far left if it is underlined. In our notation we
will enclose the distinguished element in parentheses. For example,
\[ \underline{\{\,1,6\,\}}\; \{\,8,2\,\} \; \underline{\{\,3,5\,\}} \; \{\,9,(7),4\,\}
\]
is a cyclic bipartition, which we also write
\[7,4\,\}\;\underline{\{\,1,6\,\}}\; \{\,8,2\,\} \; \underline{\{\,3,5\,\}} \;\{\,9,
\]
so that the distinguished element is equivalently the first
element in a bipartition where a block can ``wrap around." Similarly
to how we encoded signed permutations playing out on a straight table
with ordered bipartitions, we can encode the case of the circular
table with cyclic bipartitions. Here the distinguished element will
correspond to a ``distinguished guest," who sits in place number 1 at
the circular banquet table.

Let $\varphi$ be the map
encoding the circular table case. For any signed permutation $\pi$
of $[n]$, we form its image $\varphi(\pi)$ as follows (the only
difference in our algorithm is that the searches are cyclic):

\begin{enumerate}
\item[(1)] Find the least element, $\pi(i)$ (ignoring signs), that is
not already included in some subset.

\item[(2a)] If $\pi(i)$ is positive, then underline the set including
$\pi(i),$ and set $j = i+1 \mod n$.\\ While $|\pi(j)| > |\pi(j-1)|$
and $\pi(j)$ negative,\\ $\mbox{\quad}$ add $\pi(j)$ to the set
containing $\pi(i)$, and set $j = j+1 \mod n$.

\item[(2b)] If $\pi(i)$ is negative, then set $j = i-1 \mod n$.\\
While $|\pi(j)| > |\pi(j+1)|$ and $\pi(j)$ positive,\\
$\mbox{\quad}$ add $\pi(j)$ to the set containing $\pi(i)$, and set $j = j-1 \mod n$.

\item[(3)] If every element is contained in a set, then delete all
minus signs and quit. Else, go to (1).
\end{enumerate}

As an example, start with the permutation
\[ \pi = ( -7, 1, -3, 4, -2, 5, -6).
\]
The steps of the encoding are
\[ -7 \; \underline{\{\,1,-3\,\}} \; 4 \; -2 \; 5\; -6
\]
\[ -7 \; \underline{\{\,1,-3\,\}} \; \{\,4, -2\,\} \; 5\; -6
\]
and this last step, which requires us to perform a cyclic search,
\[ \underline{-7\,\}}  \; \underline{\{\,1,-3\,\}} \; \{\,4, -2\,\} \; \underline{\{\,5, -6, }
\]
Now we can drop the signs to get
\[ \underline{7\,\}} \; \underline{\{\,1,3\,\}} \; \{\,4, 2\,\} \; \underline{\{\,5, 6,}
\]
or
\[ \underline{\{\,5, 6, (7)\,\}} \; \underline{\{\,1,3\,\}} \; \{\,4, 2\,\}.
\]

Let $\varphi(\mathcal{C}_n)$ be the set of
all cyclic bipartitions corresponding to signed permutations on a
circular table. We will use this set of cyclic bipartitions along with
$\varphi_S(\mathcal{C}_n)$ to obtain the following theorem regarding the
derivative of $S(x,y,z)$.

\begin{thm}\label{thm:dSdz}
We have the following formula relating the generating function for the
straight table and the generating function for the circular table:
\begin{equation}\label{eq:dSdz}
z \frac{d}{dz} S(x,y,z) = C(x,y,z)S(x,y,z).
\end{equation}
\end{thm}

\begin{proof}
  To prove equation \eqref{eq:dSdz}, it will suffice to equate
  coefficients and prove
  \begin{equation}\label{eq:[z^n]dSdz}
    nS_n(x,y) = \sum_{i=1}^n \binom{n}{i}C_k(x,y)S_{n-i}(x,y),
  \end{equation}
  in which the polynomials
  $$
  S_n(x,y):=\sum_{\pi \in \mathcal{C}_n}\ws(\pi)\quad\text{and}\quad
  C_n(x,y):=\sum_{\pi \in \mathcal{C}_n}w(\pi)
  $$ are the coefficients of $z^n/n!$ in $S(x,y,z)$ and $C(x,y,z)$,
  respectively. We will prove \eqref{eq:[z^n]dSdz} bijectively.

  The left-hand side of equation \eqref{eq:[z^n]dSdz} can be thought of as
  counting the weights of permutations corresponding to bipartitions in
  $\varphi_S(\mathcal{C}_n)$ with a distinguished element, or a straight table with
  a distinguished guest. We simply put
  parentheses around one of the $n$ elements of the bipartition, as in
  \[
  \underline{\{\,2,8\,\} }\;\{\, 1 \,\}\;\{\, 7, (4), 3 \,\}\;
  \underline{\{\,9\,\}}\;\underline{\{\,5,6\,\}}.
  \]

  Given any bipartition with a distinguished element, we can form a pair
  $(c,s)$, where $c$ is a cyclic bipartition of a subset $A \subset [n]$
  (corresponding to one in $\varphi(\mathcal{C}_n)$), and $s$ is an ordered
  bipartition of $[n]\setminus A$ (corresponding to one in
  $\varphi_S(\mathcal{C}_n)$). We simply split the table into two pieces just
  before or after the block containing the distinguished element. If the
  block containing the distinguished element is not underlined, we make
  the split just after that block. The above example yields the pair
  \[
  \big(\; \underline{\{\,2,8\,\}}\;\{\, 1 \,\}\;\{\, 7, (4), 3
  \,\} , \, \underline{\{\,9\,\}}\;\underline{\{\,5,6\,\}} \;\big).
  \]
  If instead the block with the distinguished element is underlined, as
  in
  \[
  \{\,2,8\,\} \;\{\, 1 \,\}\;\underline{\{\, 3,(4),7
  \,\}}\;\underline{\{\,9\,\}} \;\underline{\{\,5,6\,\}},
  \]
  we make the split before the underlined block. Now the right half
  becomes the circular table, and the left half is the straight table:
  \[
  \big(\; \underline{\{\, 3,(4),7 \,\}} \;
  \underline{\{\,9\,\}}\;\underline{\{\,5,6\,\}}, \, \{\,2,8\,\}
  \;\{\, 1\, \}\;\big).
  \]
  By splitting the bipartition as we do, all of the guests retain their
  status as happy, frustrated, or napkinless. In other words, the
  product of the weights of the pair of tables equals the weight of
  original table.

  Now, given any pair $(c,s)$, where $c \in \varphi(\mathcal{C}_i)$
  and $s \in \varphi_S(\mathcal{C}_{n-i})$, we can form a
  bipartition in $\varphi_S(\mathcal{C}_n)$ with a distinguished
  element. First, in any of $\binom{n}{i}$ ways, we choose a subset $A
  = \{a_1 < a_2 < \cdots < a_i\} \subset [n]$, and replace $k$ with
  $a_k$ in $c$. For $s$, we replace $k$ with $b_k$, where
  $[n]\setminus A = \{ b_1 < b_2 < \cdots < b_{n-i} \}$. Now we
  concatenate the bipartitions. If the distinguished element of $c$ is
  in an underlined block we put the straight table on the left:
  $sc$. If the distinguished element is not in an underlined block, we
  put the straight table on the right: $cs$.
\end{proof}

Now, thanks to Corollary \ref{thm:KBLNJS} and Theorem \ref{thm:dSdz}, we have:
\[S' = \big((x-y)H\bH +q(y-1)H + p(y-1)\bH +1\big)S^2.\]
We can solve this differential equation to get the following exact formula:
\begin{equation}\label{eq:Sfunc}
S(x,y,z) = \frac{ pqy^3}{ D },
\end{equation}
where the denominator $D$ is:
\begin{equation}\label{eq:D}
\begin{aligned}
pq(y-x)e^{yz} + (qx - pqy -q^2y^2)e^{pyz} & + (px-pqy -p^2y^2)e^{qyz} \\
& + pq\big(y(y-1)^2 + x(1-yz)\big) + y^2-x.
\end{aligned}
\end{equation}
We can now obtain exact formulas for any of the other generating
functions discussed here by plugging equation \eqref{eq:Sfunc} into
the formulas of Corollary \ref{thm:KBLNJS}. In particular, we have our
main result.
\begin{thm}\label{thm:exact}
At a table for $n$ people, the probability that $i$ people are
napkinless and $j$ people are frustrated is given by the coefficient
of $x^i y^j z^k$ in the following function:
\begin{equation}\label{eq:Cfunc}
  C(x,y,z) = \frac{ pqyz\big( (x-y)e^{yz} + (qy^2 + py -x)e^{pyz} + (py^2 + qy -x)e^{qyz} + x  \big) }{ D },
\end{equation}
where $D$ is given in \eqref{eq:D}.
\end{thm}

\section{The expected number of napkinless guests and other statistics}\label{sec:expect}

With the generating function from \eqref{eq:Cfunc}, we can in
principle extract any sort of statistics related to the number of
napkinless guests or frustrated guests.  We will highlight a few statistical results that are of interest to us. In particular, we find the expected number of napkinless and frustrated guests, the variance
for each of these distributions, and the covariance for their joint
distribution.\footnote{This task is daunting by hand, but luckily we have technology to help us. Computer software such as Maple, for example, is very helpful, both with solving the differential equation leading to \eqref{eq:Sfunc}, and in obtaining residues.}

\subsection{The expected number of napkinless guests}

Recall the definition of the polynomial $C_n(x,y) = \sum_{\pi \in \mathcal{C}_n} w(\pi)$ from the proof of Theorem \ref{thm:dSdz}. Suppose we know exactly what $C_n(x,y)$ is for some
$n$. If we want to obtain the expected number of people without a
napkin, we want to compute the weighted average
\[ E_n\big(o(\pi)\big) := \sum_{\pi \in \mathcal{C}_n} p^{|\pi|_-}q^{|\pi|_+} o(\pi).
\]
In terms of the polynomial $C_n (x,1)$ (we set $y=1$ since we're not interested in frustrated diners at the moment), this just means that we differentiate with respect to $x$ and set $x =1$, since
\begin{align*}
  C'_n(x,1) &= \sum_{\pi \in \mathcal{C}_n} p^{|\pi|_-}q^{|\pi|_+} o(\pi) x^{o(\pi) -1} \\
  C'_n(1,1) &= \sum_{\pi \in \mathcal{C}_n} p^{|\pi|_-}q^{|\pi|_+} o(\pi).
\end{align*}

Therefore, we want to find the generating function for the numbers
$E_n\big(o(\pi)\big) = C'_n(1,1)$, or
\begin{align}
  E(z)
  &:= \sum_{n\geq 0} E_n\big(o(\pi)\big) z^n \nonumber\\
  &= \frac{d}{dx}\Big[ C(x,1,z) \Big]_{x=1} \nonumber\\
  &= \frac{ z\big( pq(2-z)e^z + (p^2 + pqz -1)e^{pz} + (q^2 + pqz -1)e^{qz} +1 \big) }{pq(1-z)^2}. \label{expect_gf}
\end{align}

\begin{thm}
  The expected number of napkinless guests on a circular table for $n$
  people is
 \begin{equation}\label{expect_formula}
  E_n\big(o(\pi)\big) = \frac{n}{pq}\big(1-p\exp_n(q)-q\exp_n(p)+pq\exp_n(1)\big),
 \end{equation}
  where $\exp_n(x)=\sum_{k=0}^n x^k/k!$ is the truncated exponential
  function.
\end{thm}

\begin{proof}
  Let $f(n) = n\big(1-p\exp_n(q)-q\exp_n(p)+pq\exp_n(1)\big)/(pq)$ be the right hand
  side of equation \eqref{expect_formula}, and let $F(z)$ be the ordinary
  generating function for the numbers $f(n)$. Note that
  $$nf(n+1) = (n+1)f(n) + \frac{1-p^n-q^n}{(n-1)!},
  $$
  which implies
  $$z\frac{d}{dz}\biggl[\frac{1}{z}F(z)\biggr]  =
  \frac{d}{dz}\bigl[z F(z)\bigr] + z(e^z - pe^{pz} -  qe^{qz})
  $$
  or equivalently,
  $$(1-z)F'(z) = \biggl(1+\frac{1}{z}\biggr)F(z) + z(e^z - pe^{pz} -  qe^{qz}).
  $$
  It is easy to check that $E(z)$, given by formula
  \eqref{expect_gf}, satisfies this differential equation.
\end{proof}

Formula \eqref{expect_formula} implies Sudbury's result (equation \eqref{eq:Sud}).

\begin{cor}
  The expected value of $o(\pi)$ satisfies
  $$E_n\big(o(\pi)\big) = \frac{n(1 - pe^q)(1 - qe^p)}{pq} + O\Bigl(\frac{1}{n!}\Bigr).
  $$
  In particular, when $p = q = 1/2$ we have
  $$E_n\big(o(\pi)\big) = n(2-\sqrt{e})^2 + O\Bigl(\frac{1}{n!}\Bigr).
  $$
\end{cor}
So, the answer to Winkler's problem of napkins in a random setting is
$(2-\sqrt{e})^2 \approx 0.12339675$. It was quite a bit of work for
this answer, but of course our work pays off in being able to find the
following statistics as well.

\subsection{Other statistics}
In the rest of this section we present asymptotic estimates for further statistics regarding the napkin problem. Since the formulas are less messy, we will primarily restrict our attention to the $p = q = 1/2$ case, but our approach is general.

It is straightforward to obtain asymptotic estimates for functions
with a finite number of poles such as $E(z)$. See, for example, chapter 5 of Herbert
Wilf's book \cite{Wilf}. We use the same technique to obtain all the estimates given here. We briefly outline the approach.

Suppose we have a power series $f(z) = \sum a_n z^n$, with a singularity at $z=1$ of multiplicity $m$ (and no other singularities). Then the Laurent expansion of $f$ around $1$ is: \[ f(z) = \frac{b_{-m}}{(1-z)^m} + \cdots + \frac{b_{-1}}{(1-z)} + b_0 + b_1(1-z) + b_2(1-z)^2 + \cdots. \] If we let \[ g(z) = \sum c_n z^n = \frac{b_{-m}}{(1-z)^m} + \cdots + \frac{b_{-1}}{(1-z)},\] (called the \emph{principal part} of $f$), then the function $h(z) = f(z) - g(z)$ is entire, and its coefficients vanish very quickly. Thus for large $n$, the coefficient of $z^n$ in $g(z)$ closely approximates the coefficient of $z^n$ in $f(z)$, i.e., \[ a_n \sim c_n = \binom{m+n-1}{m-1} b_{-m} + \cdots + \binom{n+1}{1}b_{-2} + b_{-1}. \] As an illustration, we apply this technique to $E(z)$ to, again, derive Sudbury's result.

We first expand $E(z)$ as a series in $u= 1-z$ to get
\begin{align*} E(u) & = \sum_{n \geq -2} b_n u^n \\
 & = \frac{ (1-u)\big( pq(1+u)e^{1-u} + (p^2 + pq(1-u) -1)e^{p(1-u)} + (q^2 + pq(1-u) -1)e^{q(1-u)} +1 \big) }{pqu^2}.
\end{align*}
We let $\overline{E}(u):= u^2 E(u)$, then set $u = 0$ to obtain \[\overline{E}(0) = b_{-2} = \frac{pqe + (p^2+pq-1)e^p + (q^2+pq-1)e^q+1}{pq}.\] Next, we differentiate $\overline{E}$ before setting $u= 0$. This gives \[ \overline{E}'(0) = b_{-1} = - \frac{pq(e + e^p + e^q) + (1+p)(p^2 + pq -1)e^p + (1+q)(q^2 + pq -1)e^q - 1}{pq}.\] Now we can obtain our estimate:
$$
a_n \sim (n + 1)b_{-2} + b_{-1} = \frac{n(1 - pe^q)(1 - qe^p)}{pq}.
$$

Now, if we want to get the variance of the distribution of napkins on the table, we need to
compute the sums of the squares of the number of napkins left on the
table, since \[ \var_n( o(\pi) ) = E_n\big( o(\pi)^2 \big) - E_n\big(
o(\pi) \big)^2,\] and we already have $E_n\big( o(\pi) \big)$. We get
\begin{align*}
  E_n\big( o(\pi)^2 \big)
  &= \sum_{\pi \in \mathcal{C}_n} p^{|\pi|_-}q^{|\pi|_+} o(\pi)^2 \\
  &= C_n''(1,1) + C_n'(1,1).
\end{align*}
Therefore the generating function for the second moment is
\[ \sum_{n\geq 0} E_n \big( o(\pi)^2 \big) z^n
= \frac{d^2}{dx^2}\Big[  C(x,1,z) \Big]_{x=1} + E(z),
\]
and all we need to find is $\displaystyle\frac{d^2}{dx^2}\Big[
C(x,1,z) \Big]_{x=1}$.  Computing, we find:
\begin{align*}
  C_{x^2}(1/2; 1,1,z)
  &= \frac{2z(2 - e^{z/2})
    \Big( e^{z/2}(1-z)^2 + (2-e^{z/2})\big( (e^{z/2}-1)^2(3-z) -e^z + 2z \big) \Big)}
  {(1-z)^3}.
\end{align*}

Using the same method, we find that the variance
is asymptotically \[ \var_n( o(\pi) ) \sim \frac{ n(1-pe^q)(1-qe^p)\big( 1 - (p^2 -pq)e^q - (q^2 - pq)e^p - pq(e + 1) \big)}{p^2q^2},\] or $n(3-e)(2-\sqrt{e})^2 \approx n(.0347631)$ for the $p=q = 1/2$ case.

To get the expectation and variance for the number of frustrated
guests, we follow the same procedure, except now we differentiate
$C(1,y,z)$ with respect to $y$. The covariance is
\[E\big(o(\pi)m(\pi)\big) - E\big(o(\pi)\big)E\big(m(\pi)\big),\] so
we differentiate with respect to $x$, then with respect to $y$ to get
the generating function for $E\big(o(\pi)m(\pi)\big)$, the only piece
we don't know. The statistics are summarized in Table
\ref{table}. Notice in particular that if the expected number of
napkinless guests is $n(2-\sqrt{e})^2$, and the expected number of
frustrated guests is $n(6\sqrt{e} - e - 7)$, then the expected number
of happy guests is $n(4-2\sqrt{e}) \approx n(.702557)$. Seventy
percent of the guests are happy!

\begin{table}[t]
\vspace{.5cm}
\begin{tabular}{|c|c|c|}
\hline
& & \\
  $X$: & $o(\pi)$ (Napkinless) & $m(\pi)$ (Frustrated) \\
\hline
& & \\
$E(X)$ & $n(2-\sqrt{e})^2\approx n(.12339675)$ & $n(6\sqrt{e} - e - 7)\approx n(.174046)$\\
& & \\
\hline
& & \\
$\var(X)$ & $n(3-e)(2-\sqrt{e})^2\approx n(.0347631)$ & $n( 6 \sqrt{e^3} - e^2 - e - 38\sqrt{e} +46)$\\
 & & $\approx n(.13138819)$ \\
\hline
& \multicolumn{2}{c|}{} \\
$\covar(X,Y)$ & \multicolumn{2}{c|}{$n\big( -(2 - \sqrt{e})( \sqrt{e^3} - 3e - 5\sqrt{e} + 12 )\big) \approx n( -.029239461)$ } \\
& \multicolumn{2}{c|}{} \\
\hline
\end{tabular}
\vspace{.5cm}
\caption{Statistics for napkinless and frustrated guests with $p = q = 1/2$.\label{table}}
\end{table}

\section{Another proof for the expected number of napkinless guests}\label{sec:alternate}

Upon reviewing a draft of this paper, Ira Gessel pointed out that the
generating function for the expected number of napkinless guests
satisfies the following differential equation (with $p=q=1/2$):
\begin{equation}\label{eq:EdN}
 E(z) = z\frac{d}{dz}\left[\frac{(2-e^{z/2})^2}{1-z} \right] = z\frac{d}{dz}\Big[ N(1,1,z)\Big],
\end{equation}
where $N(1,1,z)$ is the generating function for the proportion of
signed permutations for which neither napkin on a straight table is
taken. We will think of such permutations in terms of their image
under $\varphi_S$: specifically, those ordered bipartitions in
for which the leftmost block is underlined and the
rightmost block is not underlined. Gessel suggested that there may be
a simple combinatorial explanation for \eqref{eq:EdN}, and there is.

First, if all we want is the expected number of napkinless
guests on the circular table, then because of the symmetry of the
table, we have $E_n\big(o(\pi)\big) = np'$, where $p'$ is the
probability that any particular guest (say guest 1) has no
napkin. Upon equating coefficients, equation \eqref{eq:EdN} claims
that $|\mathcal{N}_n| = |\mathcal{C}_n^1|$, where $\mathcal{C}_n^1$ is
the set of all signed permutations for which guest 1 gets no napkin on
a circular table.

We can give a bijection between the ordered bipartitions in
$\varphi_S(\mathcal{N}_n)$ and the cyclic bipartitions in
$\varphi(\mathcal{C}_n^1)$ as follows. Given an ordered
bipartition for which the leftmost block is underlined and the
rightmost block is not underlined, there must be a guest in the middle
of the table who is napkinless. Make the leftmost such guest
distinguished, and cyclically permute the blocks of the bipartition
until the block with the distinguished guest is first or last,
depending on whether that block is underlined. For example, \[
\underline{\{\,8\,\}} \; \{\,(9),6,1\,\} \; \{\,2\,\} \; \underline{\{\,3,5\,\}}
\; \{\, 7,4\,\} \] is an ordered bipartition in $\varphi_S(\mathcal{N}_n)$
with the leftmost napkinless guest highlighted (guest 2 in this
case). We cyclically permute the blocks to get \[ \{\,2\,\} \;
\underline{\{\,3,5\,\}}\; \{\, 7,4\,\} \; \underline{\{\,8\,\}} \; \{\,(9),6,1\,\},
\] which is a cyclic bipartition where guest 1 gets no napkin on a
circular table. The inverse of this bijection is given by simply
cyclically permuting the blocks until we have an ordered bipartition
for which the distinguished guest is the leftmost guest without a
napkin, the leftmost block is underlined, and the rightmost block is
not. It is straightforward to check that we can always achieve this
state.

\section{Acknowledgements}
We would like to extend a sincere thank you to Don Knuth for bringing
this problem to our attention. His help and encouragement were
invaluable. In particular, the ideas of reducing the problem to a
straight table and keeping track of frustrated guests are his. Peter Winkler deserves credit for writing his wonderful book and for reading an early draft of this paper. We
thank the Institut Mittag-Leffler and the
organizers of the spring 2005 session in algebraic combinatorics. If
not for our time there, this work may never have been done. Lastly, we
thank Ira Gessel, whose observation motivated Section
\ref{sec:alternate}.


\begin{thebibliography}{99}

\bibitem{FoataZeil} D. Foata and D. Zeilberger, \emph{Graphical major indices}, Journal of Computational and Applied Mathematics, {\bf 68} (1996), 79--101.
\bibitem{Sloane} N. J. A. Sloane, \emph{The On-Line Encyclopedia of Integer Sequences}, (2005), published electronically at \verb"http://www.research.att.com/~njas/sequences/".
\bibitem{Sudbury} A. Sudbury, \emph{Inclusion-exclusion methods for treating annihilating and deposition processes}, Journal of Applied Probability, {\bf 39} (2002), 466--478.
\bibitem{Wilf} H. S. Wilf, \emph{generatingfunctionology}, Academic Press, San Diego, CA, 1994.
\bibitem{Winkler} P. Winkler, \emph{Mathematical Puzzles: A Connoisseur's Collection}, AK Peters, Natick, MA, 2004.

\end{thebibliography}
\end{document}